\documentclass[12pt,reqno]{amsart}

\usepackage{latexsym}
\usepackage{graphicx}

\newcommand{\C}{{\mathbb C}}

\newcommand{\Z}{{\mathbb Z}}

\renewcommand{\phi}{\varphi}

\newcommand{\ispa}[1]{\langle \,#1 \,\rangle }

\newcommand{\re}{\mathop{{\rm Re}}\nolimits}

\newcommand{\ol}{\overline}
\newcommand{\mb}{\mathbb}
\newcommand{\tr}{{\rm Tr}}

\newtheorem{theo}{{\sc Theorem}}[section]
\newtheorem{maintheo}{{\sc Theorem}}

\newtheorem{lem}[theo]{{\sc Lemma}}

\newtheorem{remark}[theo]{{\it Remark}}

\begin{document}

\title[Asymptotics of Matrix Integrals]
{Asymptotics of Matrix Integrals and Tensor Invariants of Compact
Lie Groups}

\author{Michael Stolz}
\address{Ruhr-Universit\"at Bochum, Fakult\"at f\"ur Mathematik, NA 4/32, D-44780 Bochum, Germany}
\email{michael.stolz@ruhr-uni-bochum.de}
\author{Tatsuya Tate}
\address{Graduate School of Mathematics, Nagoya University,
Furo-cho, Chikusa-ku, Nagoya, 464-8602 Japan}
\email{tate@math.nagoya-u.ac.jp}
\date{\today}

\begin{abstract}
In this paper we give an asymptotic formula for a  matrix integral
which plays a crucial role in the approach of Diaconis et al.\ to
random matrix eigenvalues. The choice of parameter for the
asymptotic analysis is motivated by an invariant theoretic
interpretation of this type of integral. For arbitrary regular
irreducible representations of arbitrary connected semisimple
compact Lie groups, we obtain an asymptotic formula for the trace
of permutation operators on the space of tensor invariants, thus
extending a result of Biane on the dimension of these spaces.

\end{abstract}

\maketitle

\section{Introduction}
\label{SIntro}

Let $G$ be a compact connected Lie group, and let $(V_{\lambda},\rho_{\lambda})$ be
an irreducible representation of $G$ with highest weight $\lambda$.
Consider the matrix integral
\begin{equation}
\label{mint} \int_{G} (\tr \rho_{\lambda}(g))^{a_{1}} (\tr
\rho_{\lambda}(g^{2}))^{a_{2}} \ldots (\tr
\rho_{\lambda}(g^{r}))^{a_{r}} ( \overline{\tr
\rho_{\lambda}(g)})^{b_{1}} \ldots ( \overline{ \tr
\rho_{\lambda}(g^{r})} )^{b_{r}} \,d\omega_G (g),
\end{equation}
where $\{a_{j}\}_{j=1}^{r}$, $\{b_{j}\}_{j=1}^{r}$ are fixed
sequences of non-negative integers and $\omega_G$ denotes
normalized Haar measure on $G$. In the case that
$(V_{\lambda},\rho_{\lambda})$ is the standard representation of
the unitary group $G = {\rm U}_n$, the integral \eqref{mint} is
nothing else than the $(a_1, \ldots, a_r, b_1, \ldots,
b_r)$-moment of the random vector
\begin{equation}
\label{vector} ({\rm Tr}(g), \ldots, {\rm Tr}(g^r), \overline{{\rm
Tr}(g)}, \ldots, \overline{{\rm Tr}(g^r)}),\end{equation} where
$g$ is chosen from ${\rm U}_n$ according to Haar measure. It has
been proven by Diaconis and Shahshahani in \cite{DS} (see also
\cite{DE}) that for $n$ large enough, \eqref{mint} coincides with
the $(a_1, \ldots, a_r, b_1, \ldots, b_r)$-moment of a vector of
independent complex Gaussian random variables of suitable
variances, and consequently, the vector \eqref{vector} converges
in distribution to this Gaussian vector as the matrix size $n$
tends to infinity. It has been observed (\cite{BR}, \cite{St})
that the Diaconis-Shahshahani result is based on the fact that the
integral \eqref{mint} can be expressed as
\begin{equation}
\label{invar}
\int_{G}\prod_{j=1}^{r}\tr(\rho_{\lambda}(g)^{j})^{a_{j}} \tr(
\rho^*_{\lambda}(g)^j)^{b_j}\,d\omega_G(g) = \tr
\left((\sigma_{k_a} \otimes \sigma_{k_b})(s,
t)|_{[V_{\lambda}^{\otimes k_a} \otimes (V_{\lambda}^*)^{\otimes
k_b}]^{G}}\right).
\end{equation}
Here we write $k_a = \sum_{j=1}^{r}ja_{j}$ and define $k_b$
analogously.  $\sigma_{k_a}$ denotes the obvious representation of
the symmetric group $\mathfrak{S}_{k_a}$ on $V_{\lambda}^{\otimes
k_a}$, and $\sigma_{k_b}$ is its analogue on
$(V_{\lambda}^*)^{\otimes k_b}$, where the representation
$(V_{\lambda}^*, \rho_{\lambda}^*)$ is contragredient to
$(V_{\lambda}, \rho_{\lambda})$. $s \in \mathfrak{S}_{k_a}$ has
cycle type $(1^{a_{1}}\ldots r^{a_{r}})$ and $t \in
\mathfrak{S}_{k_b}$ has cycle type $(1^{b_{1}}\ldots r^{b_{r}})$.
$$
 [V_{\lambda}^{\otimes k_a} \otimes (V_{\lambda}^*)^{\otimes k_b}
]^{G}= \{T \in V_{\lambda}^{\otimes k_a}  \otimes
(V_{\lambda}^*)^{\otimes k_b}: (\rho_{\lambda}^{\otimes k_a}
\otimes (\rho_{\lambda}^*)^{\otimes k_b})(g)T=T\mbox{ for all }g
\in G \} $$ is the space of invariants of the $G$-action
$\rho_{\lambda}^{\otimes k_a} \otimes (\rho_{\lambda}^*)^{\otimes
k_b}$. Diaconis and Shahshahani study the integral \eqref{mint}
for a sequence $(G_n)$ of classical groups of increasing rank,
fixing the parameters $ a = (a_1, \ldots, a_r)$ and $b = (b_1,
\ldots, b_r)$. Since they consider standard representations of
classical groups, substantial information about the right-hand
side of \eqref{invar} is available since the epoch-making work of
Weyl
(\cite{Wy}) and can be used to evaluate the left-hand side.\\

In the present paper, we depart from the Diaconis-Shahshahani
framework in two ways. Firstly, we consider arbitrary regular
irreducible representations of arbitrary compact semisimple Lie
groups. In this generality, one has {\it a priori} only poor
control of the right-hand side of \eqref{invar}. So we will
directly attack the left-hand side of this equation and thus
obtain some asymptotic information about the spaces of invariants
on the right-hand side. Secondly, we fix a group $G$ and a
representation $V_{\lambda}$, and let the parameters $k_a, k_b$ of
the tensor powers (and hence the moment parameters $a, b$) tend to
infinity. This may be thought of as a thermodynamic limit of a
particle system, rather than a random matrix limit. For our
asymptotic analysis we will use techniques which were developed by
Biane (\cite{B}), Klyachko and Kurtaran (\cite{KK}), and Tate and
Zelditch (\cite{TZ}). Our results should be compared to those of
Biane (\cite{B}) and Kuperberg (\cite{K}). An asymptotic result of
a different kind for growing tensor powers of a fixed
representation has recently been obtained by Collins and
\'{S}niady in \cite[Thm.\ 17]{CS}.\\

Here is the set-up for our main theorems: Assume that the compact
connected Lie group $G$ is semisimple, and that the highest weight
$\lambda$ of the fixed irreducible representation $(V_{\lambda},
\rho_{\lambda})$ is regular, {\it i.e.}, is in the interior of a
Weyl chamber. Fix a maximal torus $T$ in $G$. Denote by $W$ the
Weyl group. Write $\mathfrak{t}$ for the Lie algebra of $T$, and
$\mathfrak{t}^{*}$ for its dual space. $I := {\rm ker~exp~}
\subset \mathfrak{t}$ is the integral lattice, its dual $I^* := \{
\phi \in \mathfrak{t}^{*}:\ \phi(I) \subseteq \Z\}$ is the weight
lattice. Let $\Lambda^* \subset \mathfrak{t}^{*}$ be the abelian
group generated by the roots, i.e. the root lattice, and write
$\Lambda := (\Lambda^*)^*$ for its dual. It is well-known that
$\Lambda^* \subseteq I^*$, hence $I \subseteq \Lambda$, and that
the group $\Pi(G) := \Lambda / I$ is a finite abelian group. It
can be regarded as a subgroup of $T \cong \mathfrak{t} / I$. Write
$\pi$ for the canonical projection of $T$ onto $T/\Pi(G)$. For any
$\mu \in I^{*}$ write ${\rm m}_{\lambda}(\mu)$ for its
multiplicity in $V_{\lambda}$. Then, the set of all weights of
$V_{\lambda}$ is ${M}_{\lambda} := \{ \mu \in I^*: {\rm m}(\mu)
\neq 0\}$.
\\

Fix a smooth function $f > 0$ on $G$ and sequences $\alpha =
(\alpha_{1},\ldots,\alpha_{r})$ and $\beta =
(\beta_{1},\ldots,\beta_{r})$ of nonnegative integers. For a
positive integer $N$ set $a_{j} := a_{j}(N):=N\alpha_{j}$, $b_{j}
:= N\beta_{j}\ (j = 1, \ldots, r)$ and write $a = (a_1, \ldots,
a_r), b = (b_1, \ldots, b_r)$. Furthermore, we set
\begin{equation}
\label{consts}
 |\alpha| :=
\sum_{j=1}^{r}\alpha_{j},\quad k_{\alpha} := \sum_{j=1}^{r}j
\alpha_{j},\quad l_{\alpha} :=\sum_{j=1}^{r}j^{2}\alpha_{j}
\end{equation} and define $|\beta|, k_{\beta}, l_{\beta}, k_a$ etc.\
analogously. Throughout the paper we will assume that $k_{\alpha}
= k_{\beta}$.
 We
will study the integrals
\begin{eqnarray}
\label{M1} I_{N}=I_{N}(f, \alpha) &:=& \int_{G}
\prod_{j=1}^{r}(\tr (\rho_{\lambda}(g^{j})))^{N \alpha_{j}}
f(g)\,d\omega_G(g)\quad\quad\quad {\rm and}\\
\label{gen1} K_{N}=K_{N}(f, \alpha,
\beta)&:=&\int_{G}\prod_{j=1}^{r}(\tr(\rho_{\lambda}(g^{j})))^{N
\alpha_{j}}\ \overline{(\tr (\rho_{\lambda}(g^{j})))}^{N
\beta_{j}}\ f(g)\,d\omega_G(g).
\end{eqnarray}

Writing $\Phi_{+}$ for the set of positive roots of $(G,T)$, we define a polynomial $\kappa$ on $\mathfrak{t}$ by
\begin{equation}
\label{kap}
\kappa (x)=\prod_{\alpha \in \Phi_{+}}\ispa{\alpha,x},
\end{equation}
Finally, we define $A_{\lambda} \in {\rm Hom}_{\C}(\mathfrak{t}, \mathfrak{t}^*) \cong  \mathfrak{t}^*\otimes  \mathfrak{t}^* $ by
\begin{equation}
\label{mat}
A_{\lambda}=\frac{1}{\dim V_{\lambda}}
\sum_{\mu \in M_{\lambda}}m_{\lambda}(\mu)\mu \otimes \mu.
\end{equation}
Since $G$ is assumed to be semisimple, and $\lambda$ to be regular, $A_{\lambda}$ is known to be positive definite (see \cite{TZ}).\\

Now we are in a position to state our main results, to be proven
in Sections \ref{Sphase} -- \ref{Sgeneral} below.

\begin{maintheo}
\label{main1} Assume that ${\rm gcd}\{ j: \alpha_j \neq 0\} = 1.$
Then
\begin{equation}
\label{asymp1} I_{N}(f, \alpha) = \frac{(2\pi)^{d}(\dim
V_{\lambda})^{N |\alpha|} \kappa (A_{\lambda}^{-1}\rho)}{(2\pi
l_{\alpha} N)^{(\dim G)/2}\sqrt{\det A_{\lambda}}} \left( \sum_{h
\in \Pi(G)}\nu_{N k_{\alpha}\lambda}(h)f^G(h) +O(N^{-1/2})
\right),
\end{equation}
where $d$ is the number of positive roots, $\rho$ is half the sum
of the positive roots, $f^G$ is the class function $f^G(g) = \int
f(x^{-1} g x) d\omega_G(x),$ and $\nu_{N k_{\alpha} \lambda}$ is
the character on $\Pi(G)$ determined by the weight $N k_{\alpha}
\lambda$:
\begin{equation}
\label{char} \nu_{N k_{\alpha} \lambda}(h)=e^{2\pi \sqrt{-1}N
k_{\alpha} \ispa{\lambda,\psi_{h}}},
\end{equation}
$\psi_h \in \Lambda$ being a coset representative of $h \in \Lambda/I = \Pi(G)$.
\end{maintheo}

\begin{maintheo}
Suppose that ${\rm gcd}\{ j:\ \alpha_j \neq 0 {\rm ~or~} \beta_j
\neq 0\} = 1$. Then \label{main2}
\begin{equation}
K_{N}(f, \alpha, \beta) = \frac{(2\pi)^{d}(\dim V_{\lambda})^{N
(|\alpha| + |\beta|)} \kappa (A_{\lambda}^{-1}\rho)}{(2\pi N
(l_{\alpha} + l_{\beta}))^{(\dim G)/2}\sqrt{\det A_{\lambda}}}
\left( \sum_{h \in \Pi(G)}f^G(h) +O(N^{-1/2}) \right).
\end{equation}
\end{maintheo}

\begin{remark}
{\rm Note that in Theorem \ref{main1} the leading term vanishes if
$N k_{\alpha} \lambda = k_a \lambda$ is not contained in the
lattice $\Lambda^*$ and if $f \equiv 1$. In fact, in this case one
has $[V_{\lambda}^{\otimes k_a}]^G = \{ 0 \}$, hence $I_N(1,
\alpha) = 0$ by \eqref{invar}. This is because the existence of a
nontrivial invariant implies that there is a sequence
$\mu_{1},\ldots,\mu_{k_{a}}$ of weights for the irreducible
representation $V_{\lambda}$ such that $\mu_{1}+\cdots +\mu_{k_a
}=0$. Recall the well-known fact (see \cite{TZ}) that, if the
highest weight $\lambda$ is regular, then the root lattice
$\Lambda^{*}$ is spanned by the differences $\mu-\mu'$ between two
weights $\mu$, $\mu'$ for $V_{\lambda}$. This implies that
\[
1=e^{2\pi \sqrt{-1}\ispa{\mu_{1}+\cdots +\mu_{k_{a}},\psi}}
=e^{2\pi \sqrt{-1}  k_{a}\ispa{\lambda,\psi}}
\]
for any $\psi$ in the dual $\Lambda$ of the root lattice
$\Lambda^{*}$, and hence $ k_{a} \lambda \in \Lambda^{*}$.}
\end{remark}

\begin{remark}
{\rm In view of \eqref{invar}, Theorems \ref{main1} and
\ref{main2} give asymptotic formulae for the trace of permutations
on the space of tensor invariants. Specializing to the identity
permutation, one obtains the asymptotics of the dimension of these
spaces. Specifically, taking in Theorem \ref{main1} $\alpha_{j}=0$
for $j \geq 2$ and $\alpha_{1}=1$, $f \equiv 1$, and assuming that
$\lambda$ is in the root lattice $\Lambda^{*}$, we obtain an
asymptotic formula for the dimension of the space of tensor
invariants, namely
\[
\dim [V_{\lambda}^{\otimes N}]^{G}
=\frac{|\Pi (G)|(\dim V_{\lambda})^{N}\kappa (A_{\lambda}^{-1}\rho)}
{(2\pi)^{({\rm rk} G)/2}N^{(\dim G)/2}\sqrt{\det A_{\lambda}}}
(1+O(N^{-1/2})).
\]
This formula has been obtained by Biane in \cite{B} (see also
\cite{TZ}).}
\end{remark}

{\sc Acknowledgement:} Both authors have been supported by DFG via
SFB/TR 12. M.St.\ acknowledges support of JSPS in the framework of
the Japanese-German programme on infinite-dimensional harmonic
analysis. T.T.\ has been supported by JSPS Grant-in-Aid for
Scientific Research no.\ 18740089. This article was written during
his stay at Ruhr-Universit\"{a}t Bochum. He would like to thank
the people in the mathematics department, notably Prof.\ Alan
Huckleberry, for their hospitality.

\section{Matrix integrals and tensor invariants}
\label{Sinv}

Before turning to the proof of Theorems \ref{main1} and
\ref{main2}, we provide a short and self-contained proof of
equation \eqref{invar}, on which the invariant theoretic
interpretation of the integrals $I_N$ and $K_N$ is based.  Write $
H := [V_{\lambda}^{\otimes k_a} \otimes (V^*_{\lambda})^{\otimes
k_b}]^G$. Then, the orthogonal projection onto $H$ is given by
\begin{equation}
\pi_H(T) = \int_G (\rho_{\lambda}^{\otimes k_a} \otimes  (\rho_{\lambda}^*)^{\otimes k_b})(g) T\ d\omega_G(g),\quad T \in
V_{\lambda}^{\otimes k_a} \otimes (V^*_{\lambda})^{\otimes k_b}.
\end{equation}
For $A \in {\rm End}_G(V_{\lambda}^{\otimes k_a} \otimes (V^*_{\lambda})^{\otimes k_b})$, then,
$$
\tr (A|_{H})=\tr(A\pi_{H})=\tr(\pi_{H}A).
$$
We thus obtain
\begin{lem}
\label{tr11} For any
$A \in {\rm End}_G(V_{\lambda}^{\otimes k_a} \otimes (V^*_{\lambda})^{\otimes k_b})$, we have
$$
\tr (A|_{H})= \int_G \tr\left((\rho_{\lambda}^{\otimes k_a}
\otimes  (\rho_{\lambda}^*)^{\otimes k_b})(g)\ A\right)
d\omega_G(g) = \int_G \tr\left( A\ (\rho_{\lambda}^{\otimes k_a}
\otimes  (\rho_{\lambda}^*)^{\otimes k_b})(g)\right)
d\omega_G(g).$$
\end{lem}

In view of this, \eqref{invar} is implied by the following well-known lemma (see \cite{DS}, \cite{Ra}):

\begin{lem} Let $V$ be a $d$-dimensional complex vector space, $B \in {\rm End}_{\C}(V)$, and  $s \in \mathfrak{S}_{k}$ a permutation of type
$(1^{a_1} 2^{a_2}\ldots r^{a_r})$, hence $k = k_a = \sum_{j=1}^r j a_j$. Let $\mathfrak{S}_{k}$ act on $V^{\otimes k}$ via
$\sigma_k(s)(  \otimes_{i = 1}^k v_i ) := \otimes_{i=1}^k v_{i s^{-1}}.$
Then the trace of $B^{\otimes k} \sigma_k(s) \in {\rm End}_{\C}(V^{\otimes k})$ is
$$ \tr\left( B^{\otimes k} \sigma_k(s) \right) = \prod_{j=1}^r \tr(B^j)^{a_j}.$$
\end{lem}

\begin{proof}
Fix an inner product $\ispa{,}$ on $V$.
This induces an inner product on $V^{\otimes k}$,
which is also denoted by $\ispa{,}$.
Let $e_{1},\ldots,e_{d}$
be an orthonormal basis for $V$. Write ${\mathcal F}$  for the set of maps from $\{ 1, \ldots, k\}$ to $ \{ 1, \ldots, d\}$,
${\mathcal F}_S$ for the restrictions to a subset $S$ of $\{ 1, \ldots, k\}$, $e_{\phi} := \otimes_{i=1}^k e_{\phi(i)}.$
Then $\{ e_{\phi}: \phi \in {\mathcal F}\}$ is an orthonormal basis of $V^{\otimes k}$. Write
$$ s = \prod_{j=1}^r \prod_{i=1}^{a_j} \zeta_i^j,$$ where $ \{ \zeta_i^j:\ i = 1, \ldots, a_j\}$ are the cycles of length $j$ in $s$.
 Furthermore, for any $t \in \mathfrak{S}_{k}$,
write $[t] := \{ \nu = 1, \ldots, k: \nu t \neq  \nu\}.$
Then
\begin{eqnarray*}
&&\tr(B^{\otimes k} \sigma_k(s) ) = \sum_{\phi \in {\mathcal F}} \left\langle B^{\otimes k} \sigma_k(s) e_{\phi} , e_{\phi} \right\rangle =
 \sum_{\phi \in {\mathcal F}} \prod_{j = 1}^k \left\langle B e_{j s^{-1} \phi}, e_{j \phi}\right\rangle
\\ &=& \sum_{\phi \in {\mathcal F}} \prod_{j=1}^r \prod_{i=1}^{a_j} \prod_{l \in [ \zeta_i^j]} \left\langle B e_{l (\zeta_i^j)^{-1}\phi}, e_{l\phi}\right\rangle
=  \prod_{j=1}^r \prod_{i=1}^{a_j} \sum_{\phi \in {\mathcal F}_{\![\zeta_i^j]}}  \prod_{l \in [ \zeta_i^j]} \left\langle B e_{l (\zeta_i^j)^{-1}\phi}, e_{l\phi}
\right\rangle\\
&=& \prod_{j=1}^r \prod_{i=1}^{a_j} \tr(B^j) = \prod_{j=1}^r \tr(B^j)^{a_j}.
\end{eqnarray*}

\end{proof}

\section{Phase function for the matrix integral}
\label{Sphase}

We start by rewriting the integral \eqref{M1} using Weyl's
integration formula, assuming for simplicity that $f$ is a class
function:
\begin{equation}
\label{weyli} I_{N}(f, \alpha) =\frac{1}{|W|}
\int_{T}\prod_{j=1}^{r}\tr(\rho_{\lambda}(t^{j}))^{N
\alpha_{j}}f(t)|\Delta (t)|^{2}\,dt,
\end{equation}
where $T \subset G$ is a maximal torus, $dt$ is Haar measure on
$T$, normalized as a probability measure and $\Delta(t)$ is the
Weyl denominator. We define the following function on the
complexified Lie algebra $\mathfrak{t}^{\mb{C}}$:
\begin{equation}
\label{char1}
k(w)=\sum_{\mu \in M_{\lambda}}m_{\lambda}(\mu)e^{2\pi\ispa{\mu,w}},
\end{equation}
where $w=\tau +\sqrt{-1}\varphi \in \mathfrak{t}
\oplus \sqrt{-1}\mathfrak{t}=\mathfrak{t}^{\mb{C}}$,
and linear forms in $\mathfrak{t}^{*}$ are
extended complex linearly to $\mathfrak{t}^{\mb{C}}$.
Note that the restriction of $k$ to $\sqrt{-1}\mathfrak{t}$ is
essentially the character of $V_{\lambda}$.

Let $d\varphi$ denote Lebesgue measure on $\mathfrak{t}$,
normalized so that the fundamental domain $T_{o}$ of the integral
lattice $I$ has volume $1$. Then the integral $I_{N}$ can be
written in the form:
\begin{equation}
\label{alg1} I_{N}=\frac{1}{|W|}\int_{T_{o}}
F(\sqrt{-1}\varphi)^{N} f(\varphi)|\Delta(\varphi)|^{2}\,d\varphi,
\quad
F(\sqrt{-1}\varphi)=\prod_{j=1}^{r}k(\sqrt{-1}j\varphi)^{\alpha_{j}},
\end{equation}
where the class function $f$ is, through the exponential map,
regarded as a function on $\mathfrak{t}$.

\begin{lem}
\label{kernel1}
We have the inequality
\begin{equation}
\label{realp} \left| F(\sqrt{-1}\varphi) \right| \leq
k(0)^{|\alpha|}=(\dim V_{\lambda})^{|\alpha|}.
\end{equation}
Equality holds if and only if $\varphi$ is in the dual lattice
$\Lambda$ of the root lattice $\Lambda^{*}$.
\end{lem}
\begin{proof} Similar to the proof of Lemma 1.4 in \cite{TZ}.
Note that the assumption about the gcd is used here.\end{proof}

By Lemma \ref{kernel1}, the integral \eqref{weyli} or \eqref{alg1}
is localized on $\ker \pi = \Pi(G) = \Lambda/I$. Let $g$ be a
smooth cut-off function on $T$ around the unit such that the
support of $g$ does not contain any element in $\ker \pi$ other
than the unit. Then, the translate $g_{h}(t)=g(h^{-1} t)$ is a
cut-off function around $h \in \ker \pi$. Let $\varphi \in
\Lambda$. Then it is easy to see that
\begin{equation}
\label{equi1}
k(\sqrt{-1}j\varphi)=e^{2\pi \sqrt{-1}j\ispa{\lambda,\varphi}}k(0),
\end{equation}
which is not zero. Thus, around each $\varphi \in \Lambda$, we can take a
branch of the logarithm to define the following function $\Phi$:
\begin{equation}
\label{phase1} \Phi(w):= \sum_{j=1}^{r}\alpha_{j}\log k(jw),\quad
w=\tau +\sqrt{-1}\varphi \in \mathfrak{t}^{\mb{C}},
\end{equation}
where $\varphi$ varies in a neighborhood of a point in $\Lambda$.
Then,  by Lemma \ref{kernel1}, we can write the integral $I_{N}$
as follows:
\begin{equation}
\label{sumi} I_{N}=\frac{1}{|W|} e^{N \Phi(0)} \sum_{h \in \ker
\pi} \int_{\mathfrak{t}} e^{N( \Phi(\sqrt{-1}\varphi) -
\Phi(0))}g_{h}(\varphi)f(\varphi) |\Delta (\varphi)|^{2}\,d\varphi
\end{equation}
plus a term of order $O(e^{-cN})$ for some $c > 0$. To compute
each of the integrals in the sum in \eqref{sumi}, we note that
\[
\Delta (\varphi + \psi_h)=e^{2\pi
\sqrt{-1}\ispa{\rho,\psi_{h}}}\Delta(\varphi)
\]
for each $h \in \ker \pi \cong \Lambda /I$, where $\psi_{h} \in
\Lambda$ satisfies $\exp (\psi_{h})=h$ and $\rho$ is half the sum
of the positive roots. Furthermore, by \eqref{equi1}, we have
\[
\Phi(\sqrt{-1}(\varphi + \psi_{h}))=\Phi(\sqrt{-1}\varphi) + 2\pi
\sqrt{-1}\ k_{\alpha} \ispa{\lambda,\psi_{h}}.
\]
Therefore, we obtain:
\begin{lem}
\label{decomp1}
We have
\begin{equation}
\label{dec1} I_{N}=\frac{1}{|W|}\sum_{h \in \ker \pi} \nu_{N
k_{\alpha} \lambda}(h) \int e^{N\Phi
(\sqrt{-1}\varphi)}g(\varphi)f_{h}(\varphi)|\Delta
(\varphi)|^{2}\,d\varphi
\end{equation}
plus a term of order $O(e^{-cN})$, where $g(\varphi)$ is a cut-off
function around $\varphi=0$, and $f_{h}(\varphi)=f(\varphi +
\psi_{h})$ with a representative $\psi_{h} \in \Lambda$ for $h \in
\ker \pi$.
\end{lem}

Next, we compute the first and second derivatives of the phase
function $\Phi$ at points in $\ker \pi$.
\begin{lem}
\label{critic1} Any $\varphi \in \Lambda$ is a critical point of
$\Phi$. Furthermore, the negative of the Hessian of $\Phi$,
$H(\varphi):\mathfrak{t} \to \mathfrak{t}^{*}$, is given by
$H(\varphi)=(2\pi)^{2}\ l_{\alpha} A_{\lambda}$ which is
independent of $\varphi \in \Lambda$ and is positive definite,
where $l_{\alpha}$ and $A_{\lambda}$ are defined in \eqref{consts}
and \eqref{mat}, respectively.
\end{lem}
\begin{proof}
That each $\varphi \in \Lambda$ is a critical point of $\Phi$ is
proven by the fact that
\[
\sum_{\mu \in M_{\lambda}}m_{\lambda}(\mu)\mu=0
\]
because the left-hand side of the above is a $W$-invariant vector
and $G$ is assumed to be semi-simple. The fact that the linear map
$A_{\lambda}:\mathfrak{t} \to \mathfrak{t}^{*}$ is positive
definite is proven in \cite{TZ}.\end{proof}

Therefore, what we need to do is to find asymptotics of each of
the integrals in \eqref{dec1}, each of which is an integral of
functions supported around the origin. However, there is a
difficulty: for each integral, the origin is the unique critical
point for the complex phase function $\Phi$, but the Weyl
denominator $\Delta$ vanishes at the origin. In the next section,
we will use ideas of Biane (\cite{B}) and Klyachko-Kurtaran
(\cite{KK}) to circumvent this problem.

\section{The method of Biane and Klyachko-Kurtaran}
\label{Sbiane}

In this section we study the following integral, which is a slight
generalization of the integral on the right-hand side of
\eqref{dec1}:
\begin{equation}
\label{intG}
J_{N}:=\int_{T}
e^{N\Phi (t)}g(t)|\Delta(t)|^{2}\,dt,
\end{equation}
where $dt$ is normalized Haar measure on the maximal torus $T$ and
$\Delta$ denotes the Weyl denominator. $g$ is a compactly
supported smooth function which does not vanish in the unit
element of $T$.
 We fix a $W$-invariant inner product on
$\mathfrak{t}$ such that the volume of the parallelotope
determined by an orthonormal basis equals $1$. We assume that the
complex valued smooth phase function $\Phi$, which we view as a
function on the complexified Lie algebra $\mathfrak{t}^{\mb{C}}$,
satisfies the following conditions:

\begin{itemize}
\item[(i)] $\re (\Phi(\sqrt{-1}\varphi)) \leq \Phi (0)$,  with
equality if and only if $\phi =0$; \item[(ii)] the origin is a
critical point of $\Phi$; \item[(iii)] $H=-\partial^{2}\Phi
(0):\mathfrak{t} \to \mathfrak{t}^{*}$, the negative of the
Hessian of $\Phi$ at the origin, is positive definite; \item[(iv)]
in the Taylor expansion
\begin{equation} \label{taylorP} \Phi
(\sqrt{-1}\varphi)-\Phi(0)=-\ispa{H\varphi,\,\varphi}/2
-\sqrt{-1}\ \Theta(\varphi)+R_{4}(\varphi)
\end{equation} up to fourth order,
 $\Theta$ is a real valued homogenous
polynomial of degree 3 \item[(v)] the linear map $H$ commutes with
the action of the Weyl group $W$.
\end{itemize}
Note that the phase functions of the integrals $I_N$ and $K_N$
(see Sections \ref{Sphase} and \ref{Sgeneral}) satisfy these
conditions. We are now in a position to state the main result of
this section:

\begin{theo}
\label{vanish} Under the above conditions on $\Phi$ and $g$ one
has
\begin{equation}
J_{N}= \left( \frac{2\pi} {N} \right)^{(\dim G)/2}\
\frac{(2\pi)^{d}g(0)e^{N\Phi(0)}|W|} {\sqrt{\det H}}\ \kappa
(H^{-1}\rho)\ (1+O(N^{-1/2})),
\end{equation}
where $|W|$ is the order of the Weyl group $W$,
$d$ is the number of positive roots, $\rho$ is half the sum of the positive roots,
and the polynomial $\kappa$ on $\mathfrak{t}$
is defined in \eqref{kap}.
\end{theo}

\begin{proof}
First of all, we normalize Lebesgue measure on $\mathfrak{t}$ so
that the volume of the fundamental domain $T_{o}$ of the lattice
$I$ equals $1$, and write
\begin{equation}
\label{aux01}
J_{N}=e^{N\Phi(0)}\int_{T_{o}}e^{N[\Phi(\sqrt{-1}\varphi)-\Phi(0)]}
g(\varphi)|\Delta(\varphi)|^{2}\,d\varphi.
\end{equation}
We may regard $g$ as a function on $\mathfrak{t}$ with arbitrarily
small compact support around the origin, since by (i) the
integrand in \eqref{aux01} is bounded by $e^{-cN}$ (with a
constant $c>0$) outside a compact neighborhood of the origin. (i)
and (iv) imply that for fixed $a>0$ one can choose $b>0$ such that
$\re \Phi(\sqrt{-1}\varphi) -\Phi(0) \leq
-b\ispa{H\varphi,\,\varphi}$ for $|\varphi| \leq a$. Substituting
\eqref{taylorP} into $J_{N}$ and changing the variable $\varphi$
to $N^{- 1/2}\varphi$, we have
\begin{equation}
\label{aux11}
\begin{split}
J_{N}=& N^{- ({\rm rk} G)/2}e^{N\Phi(0)} \\
&\times \int e^{-\ispa{H\varphi,\,\varphi}/2\ -\sqrt{-1}N\ \Theta
(N^{-1/2}\varphi)\ +N R_{4}(N^{-1/2}\varphi)}
g(N^{-1/2}\varphi)|\Delta(N^{-1/2}\varphi)|^{2}\,d\varphi.
\end{split}
\end{equation}
As in \cite{B}, \cite{TZ}, using the identity $\Delta
=\prod_{\alpha \in \Phi_+} (e^{\pi \sqrt{-1}\alpha}-e^{-\pi
\sqrt{-1}\alpha})$, it is easy to see that
\[
\Delta(N^{-1/2}\varphi)=(2\pi \sqrt{-1})^{d}N^{-d/2}\kappa(\varphi)
(1+O(N^{-1}|\varphi|^{2})).
\]
We note that $g(N^{-1/2}\varphi)=g(0)(1+O(N^{-1/2}|\varphi|))$.
Substituting these formulas into \eqref{aux11}  and introducing a
smooth cut-off function $\chi$  such that $\chi = 1$ around the
support of $g$, we obtain
\begin{equation}
\label{aux12}
\begin{split}
J_{N}=&\frac{(2\pi)^{2d}g(0)e^{N\Phi(0)}}{N^{d+({\rm rk} G)/2}} \\
& \times \int e^{-\ispa{H\varphi,\,\varphi}/2- \sqrt{-1} N^{-1/2}
\Theta(\varphi)+NR_{4}(N^{-1/2}\varphi)}
\chi(N^{-1/2}\varphi)|\kappa(\varphi)|^{2}\,d\varphi \
(1+O(N^{-1/2})).
\end{split}
\end{equation}
Here, we note that $g$ and its derivatives are bounded on
$\mathfrak{t}$, and that the exponential in the integrand is
bounded by $e^{-b\ispa{H\varphi, \varphi}}$ if $|\varphi|/N^{1/2}
\leq a$. As indicated above, we may assume that
$g(N^{-1/2}\varphi)=0$ for $|\varphi|/N^{1/2} \geq a$.

Next, as in \cite{TZ}, we divide the integral in \eqref{aux12}
into several parts as follows. We set
$E_{N}(\varphi):=e^{-\ispa{H\varphi,\varphi}/2- \sqrt{-1} N^{-1/2}
\Theta(\varphi)}$ and write
\begin{equation}
J_{N}=\frac{(2\pi)^{2d}g(0)e^{N\Phi(0)}}{N^{d+ ({\rm rk} G)/2}}
\left( \sum_{j=1}^{4}I_{j}(N) \right)(1+O(N^{-1/2})),
\end{equation}
where the integrals $I_{j}(N)$, $j=1,2,3,4$ are given by
\begin{equation}
\begin{gathered}
I_{1}(N)=\int E_{N}(\varphi)
|\kappa (\varphi)|^{2}\,d\varphi,\\
I_{2}(N)=\int E_{N}(\varphi)
(e^{NR_{4}(N^{-1/2}\varphi)}-1)\chi(N^{-1/4}\varphi)
|\kappa(\varphi)|^{2}\,d\varphi,\\
I_{3}(N)=\int E_{N}(\varphi)
e^{NR_{4}(N^{-1/2}\varphi)}
(1-\chi(N^{-1/4}\varphi))\chi(N^{-1/2}\varphi)|\kappa (\varphi)|^{2}\,d\varphi,\\
I_{4}(N)=\int E_{N}(\varphi)
(\chi(N^{-1/4}\varphi)-1)|\kappa (\varphi)|^{2}\,d\varphi,
\end{gathered}
\end{equation}
where we used the relation
$\chi(N^{-1/4}\varphi)\chi(N^{-1/2}\varphi)=\chi(N^{-1/4}\varphi)$ for
sufficiently large $N$.

Now, the integrand in $I_{2}(N)$ vanishes for $|\varphi| \geq
cN^{1/4}$, and thus we have $e^{NR_{4}(N^{-1/2}\varphi)}=O(1)$ and
$NR_{4}(N^{-1/2}\varphi)=|\varphi|^{4}O(1/N)$. Hence, $(e^{N
R_4(N^{-1/2}\phi)} - 1)$ is bounded by $c |\phi|^4
O(\frac{1}{N}).$ But we have $|E_{N}(\varphi)| =
e^{-\ispa{H\varphi,\varphi}/2}$, and hence
$|\varphi|^{4}E_{N}(\varphi)$ is integrable uniformly in $N$. Thus
we have $I_{2}(N)=O(1/N)$. As to $I_{3}(N)$, the function
$E_{N}(\varphi)e^{NR_{4}(N^{-1/2}\varphi)}$ is dominated by
$e^{-b\ispa{H\varphi,\,\varphi}}$ wherever the integrand does not
vanish. Since $\chi(N^{-1/4}\varphi)=1$ for $|\varphi| \leq c
N^{1/4}$ for some $c>0$, we easily have
$I_{3}(N)=O(N^{c_{1}}e^{-c_{2}N^{1/2}})$ for some $c_{1},c_{2}>0$.
Similarly, we have $I_{4}(N)=O(N^{c_{1}}e^{-c_{2}N^{1/2}})$.

Finally, we consider the integral $I_{1}(N)$. Note that $e^{-
\sqrt{-1}\ \Theta(\varphi) N^{-1/2}}=1+O(|\varphi|^{3}/N^{1/2})$
since $\Theta(\varphi)$ is real. Thus, invoking the identity
$$\int_{\mathfrak{t}}
e^{-\ispa{Hx,x}/2}|\kappa(x)|^{2}\,dx =\frac{(2\pi)^{({\rm rk}
G)/2}\ |W|\ \kappa(H^{-1}\rho)} {\sqrt{\det H}}, $$ known as
Mehta's conjecture and proven in \cite{Mc} and \cite{Op}, we have
\begin{equation}
\label{aux31} I_{1}(N)=\int e^{-\ispa{H\varphi,\,\varphi}}|\kappa
(\varphi)|^{2}\,d\varphi (1+O(N^{-1/2})) = \frac{(2\pi)^{({\rm rk}
G)/2}\ |W|\ \kappa (H^{-1}\rho)}{\sqrt{\det H}}(1+O(N^{-1/2})),
\end{equation}
which completes the proof.
\end{proof}

\section{Proof of the main results}
\label{Sgeneral}

Theorem \ref{main1} is a direct application of Theorem
\ref{vanish} with Lemma \ref{critic1} to the integral on the
right-hand side of \eqref{dec1}. To prove Theorem \ref{main2}, we
proceed as in Section \ref{Sphase} and use the Weyl integration
formula to obtain
\begin{equation}
\label{gen2} K_{N}=\frac{1}{|W|}\int_{T} J(t)^{N}
f(t)|\Delta(t)|^{2}\,dt,\quad J(t)=\prod_{j=1}^{r} \tr
(\rho_{\lambda}(t^{j}))^{\alpha_{j}}
\overline{\tr(\rho_{\lambda}(t^{j}))}^{\beta_{j}}.
\end{equation}

We also apply Theorem \ref{vanish} to find asymptotics of
the integral $K_{N}$ as follows.
As in Lemma \ref{kernel1}, we have
\begin{equation}
\label{kernel2} |J(\sqrt{-1}\varphi)| \leq k(0)^{|\alpha| +
|\beta|}=(\dim V_{\lambda})^{|\alpha| + |\beta|},
\end{equation}
with equality if and only if $\varphi \in \Lambda$. We define
\begin{equation}
\label{Gphase}
\Psi (w)=\sum_{j=1}^{s}[\alpha_{j}\log k(jw)+\beta_{j} \log k(j\ol{w})],\quad
w=\tau +\sqrt{-1}\varphi \in \mathfrak{t}^{\mb{C}},
\end{equation}
around each $\varphi \in \Lambda$.
Next, we need to compute the Hessian of $\Psi$ at $\varphi \in \Lambda$.
We note that, for $\varphi \in \Lambda$, we have
\[
\begin{gathered}
k(\sqrt{-1}j\varphi)=e^{2\pi \sqrt{-1}j\ispa{\lambda,\varphi}}k(0), \quad
(\partial k)(\sqrt{-1}j\varphi) =0,\\
\partial^{2}(\log k(\sqrt{-1}j\varphi))=-4\pi^{2}j^{2}A_{\lambda}.
\end{gathered}
\]
Therefore, we have the following lemma:
\begin{lem}
Each $\varphi \in \Lambda$ is a critical point of $\Psi$. The
negative of the Hessian of $\Psi$, denoted $D(\varphi)$, is given
by $D(\varphi)=(2\pi)^{2} (l_{\alpha} + l_{\beta}) A_{\lambda}$,
and hence $D=D(\varphi)$ does not depend on $\varphi \in \Lambda$,
and is positive definite.
\end{lem}

This time, the term involving $2\pi \sqrt{-1} j \ispa{\lambda,\psi_{h}}$ disappears
because of the assumption $\sum j\alpha_{j}=\sum j\beta_{j}$.
We thus have
\[
K_{N}=\frac{1}{|W|}\sum_{h \in \ker \pi}\int_{\mathfrak{t}}
e^{N\Psi (\sqrt{-1}\varphi)}g(\varphi)
f_{h}(\varphi)|\Delta (\varphi)|^{2}\,d\varphi,
\]
where $f_{h}$ is defined in Lemma \ref{decomp1}.
The assumption of Theorem \ref{vanish} is satisfied by the phase function $\Psi$,
and applying it, we obtain Theorem \ref{main2}.

\end{document}